\documentclass[12pt]{article}
\usepackage{amsmath,amssymb}
\usepackage{amsthm}
\usepackage{array}
\newtheorem{thm}{Theorem}[section]
\newtheorem{prop}[thm]{Proposition}
\newtheorem{cor}[thm]{Corollary}
\newtheorem{lemma}[thm]{Lemma}

\theoremstyle{note}
\newtheorem{defn}[thm]{Definition}

\newtheorem{exam}[thm]{Example}

\def\a{\mathfrak a}
\def\b{\mathfrak b}
\def\A{\mathfrak A}

\def\K{\mathfrak K}
\def\E{\mathfrak E}
\def\SS{\mathfrak S}
\def\M{\mathfrak M}
\def\H{\mathfrak H}
\def\N{\mathfrak N}
\def\rit#1{{\mbox{\rm #1}}}
\def\modx#1#2{\equiv#1\hspace{-1mm}\mod #2}
\def\itemx#1{\item[{\rm(#1)}]}

\begin{document}
\title{Singular values of some modular functions}
\maketitle
\begin{center}
 N{\sc oburo} I{\sc shii} and M{\sc aho} K{\sc obayashi} \end{center}
\section{Introduction}
 For a positive integer $N$, let $\Gamma_0(N)$ and $\Gamma_1(N)$ be the subgroups of $\rit{SL}_2(\mathbf Z)$ defined by
\[
\begin{split}
\Gamma_0(N)&=\left\{\left. \begin{pmatrix} a & b \\ c & d \end{pmatrix}\in \rit{SL}_2(\mathbf Z)~\right |~ c \equiv 0 \mod N \right\},\\
\Gamma_1(N)&=\left\{\left. \begin{pmatrix} a & b \\ c & d \end{pmatrix}\in \rit{SL}_2(\mathbf Z)~\right |~ a-1\equiv c \equiv 0 \mod N \right\}.
\end{split}
\]
We denote by $A_1(N)$ and $A_0(N)$ the modular function fields with respect to $\Gamma_1(N)$ and $\Gamma_0(N)$ respectively. Let $\mathfrak E$ be a set of triples of integers $\mathfrak a=[a_1,a_2,a_3]$ with the properties $0<a_i\leq N/2$ and $a_i\ne a_j$ for $i\ne j$. For an element $\tau$ of complex upper half plane $\mathfrak H$, we denote by $L_\tau$ the lattice in $\mathbf C$ generated by $1$ and $\tau$. Let $\wp(z;L_\tau)$ be the Weierstrass $\wp$-function relative to the lattice $L_\tau$. For $\a\in\mathfrak E$, consider a function $W_{\a}(\tau)$ on $\mathfrak H$ defined by
 \[
W_{\a}(\tau)=\frac{\wp (a_1/N;\tau)-\wp (a_3/N;\tau)}{\wp (a_2/N;\tau)-\wp (a_3/N;\tau)}.
\]
This function is a modular function with respect to $\Gamma_1(N)$, referred in Chapter 18, \S6 of Lang \cite{LA}. He pointed out that it is interesting to investigate its special values at imaginary quadratic points. In \cite{II} and \cite{I1}, to construct generators of $A_1(N)$ and $A_0(N)$, we used the function $W_\a(\tau)$ and the function $T_{\a_1,\a_2}(\tau)$ which is the trace of the product $W_{\a_1}W_{\a_2}$ ($\a_i\in\E$) relative to the extension $A_1(N)/A_0(N)$. Further we provided an explicit representation of the modular $j$-function $j(\tau)$ with those generators. In this article, we study the properties of singular values of $W_\a$ and those of a function $T_{\A,F}$ which is a generalization of the function $T_{\a_1,\a_2}$. See \S 2 for the precise definition of $T_{\A,F}$. Our results in this article are as follows. In Theorem~\ref{singv} and Corollary \ref{cor1} we prove, for imaginary quadratic points $\alpha\in\H$ and sets $\a,\A$ satisfying some conditions, that singular values $W_\a(\alpha)$ are units of the ray class field $\K_N$ modulo $N$ over $K$ and that singular values $T_{\A,F}(\alpha)$ are algebraic integers in $\K_N$. In particular, consider the triples $\a_1=[2,3,1]$ and $\a_2=[2,5,1]$. Then we prove in Theorem \ref{generator} that $W_{\a_1}(\alpha)$ and $W_{\a_2}(\alpha)$ generate $\K_N$ over the field $K(\exp(2\pi i/N))$. Let $A_0(N)_{\mathbf Q}$ be the subfield of $A_0(N)$ consisting of modular functions with Fourier coefficients in $\mathbf Q$. In Proposition \ref{progen1} we show for prime numbers $N$ that $A_0(N)_{\mathbf Q}=\mathbf Q(T_{\a_1},T_{\a_2})$ $=\mathbf Q(T_{\a_i},T_{\a_1,\a_2})$ $(i=1,2)$. Further put $\A_0=[\a_1,\a_2]$ and $F_0=X_1^mX_2^n$ for non-negative integers $m$ and $n$. In Theorem \ref{progen3}, without the assumption $N$ are prime, we show that $A_0(N)_{\mathbf Q}=\mathbf Q(j,T_{\A_0,F_0})$. We deduce from those results that singular values of those functions generate ring class fields over $K$ (see Theorem \ref{singt}). Finally in \S5 we study class polynomials of $T_{\A,F}$ with respect to Schertz $N$-systems. \par
 In the followings, for a function $f(\tau)$ and a matrix $A=\begin{pmatrix}a&b\\c&d\end{pmatrix}\in\rit{SL}_2(\mathbf Z)$, we shall denote
\[
f[A]_2=f\left(\frac{a\tau+b}{c\tau+d}\right)(c\tau+d)^{-2}\text{ and }f\circ A=f\left(\frac{a\tau+b}{c\tau+d}\right).
\]
\section{Modular functions $W_\a(\tau)$ and $T_{\A,F}(\tau)$}
Let $W_{\a}(\tau)$ be the function defined in \S 1. In \cite{II}, we showed the function $W_{\a}$ is a modular function with respect to $\Gamma_1(N)$ and it has neither zeros nor poles on $\mathfrak H$.  Let us consider the factor group $G(N)=\Gamma_0(N)/\{\pm E_2\}\Gamma_1(N)$, where $E_2$ is the unit matrix. Put $\SS_N= (\mathbf Z/N\mathbf Z)^\times/\{\pm 1\}$. Then $$G(N)\cong\left\{\left.\begin{pmatrix}\lambda^{-1} & 0\\0&\lambda\end{pmatrix} \right | \lambda\in\SS_N\right\}.$$ For $\lambda \in \SS_N$, let $M_\lambda\in\Gamma_0(N)$ such that $M_\lambda\equiv \big(\begin{smallmatrix}\lambda^{-1}&0\\0&\lambda\end{smallmatrix}\big)\mod N$.  For a tuple $\A=[\mathfrak a_1,\dots,\mathfrak a_n]~(\a_i\in\E)$ and a polynomial $F=F(X_1,X_2,\dots,X_n)$ $\in\mathbf Z[X_1,X_2,\dots,X_n]$, we define a function
\[
T_{\A,F}(\tau)=\sum_{\lambda\in\SS_N}F(W_{\a_1}\circ M_\lambda,\cdots, W_{\a_n}\circ M_\lambda).
\]
Then obviously $T_{\A,F}(\tau)$ is a modular function with respect to $\Gamma_0(N)$ and has no poles on $\H$. For $\lambda \in\SS_N,\a=[a_1,a_2,a_3]\in\E$, define an element $\lambda\a$ of $\E$ by
 $$\lambda\a=[\{\lambda a_1\},\{\lambda a_2\},\{\lambda a_3\}],$$ where $\{\lambda a_i\}$ is the integer such that $\{\lambda a_i\}\equiv\pm\lambda a_i \mod N,~0<\{\lambda a_i\}\leq\frac N 2$. 
\begin{prop}\label{pro1}
\begin{enumerate}
\itemx i $W_\a(M_\lambda\tau)=W_{\lambda\a}(\tau)$.\\
\itemx {ii}~$\displaystyle T_{\A,F}(\tau)=\sum_{\lambda\in\SS_N}F(W_{\lambda\a_1}(\tau),\cdots,W_{\lambda\a_n}(\tau))$.

\end{enumerate}
\end{prop}
\begin{proof}
The assertion (i) is showed in \S2 of \cite{II}. The assertion (ii) is obvious from (i).
\end{proof}
 We denote by $T_\a$ and $T_{\a_1,\a_2}$ the function $T_{\A,F}$ with $\A=[\a],F=X_1$ and $\A=[\a_1,\a_2],F=X_1X_2$ respectively. 
\section{Modular equations}
Let $j$ be the modular $j$-function. Let $\Gamma$ be a subgroup of $\rit{SL}_2(\mathbb Z)$ of finite index. For a modular function $f$ with respect to $\Gamma$, we define the modular equation of $f$ relative to $j$ by
\[
\Phi[f](X,j)=\prod_B(X-f\circ B),\]
where $B$ runs over a transversal of the coset decomposition of $\rit{SL}_2(\mathbf Z)$ by $\Gamma$. Obviously the coefficients of $\Phi[f](X,j)$ with respect to $X$ are in $\mathbf C(j)$. If $f$ has no poles on $\H$, then the coefficients of $\Phi[f](X,j)$ are polynomials of $j$. Hereafter to avoid tedious notation, we denote by $\Phi_{\A,F}(X,j)$ the equation $\Phi[T_{\A,F}](X,j)$. Since $W_\a$ and $T_{\A,F}$ have no poles on $\H$, we have $\Phi[W_\a](X,j),~\Phi_{\A,F}(X,j)\in\mathbf C[j][X]$. We shall show that $\Phi[W_\a](X,j)$ and $\Phi_{\A,F}(X,j)\in\mathbf Z[j][X]$ under some conditions imposed on $N$ and $\A$. For a positive divisor $t$ of $N$, let $\Theta_t$ be a set of $\varphi((t,N/t))$ pairs of integers $(u,v)$ such that
$(u,t)=1,~uv\equiv 1\mod t$ and u are inequivalent to each other modulo $(t,N/t)$. For $(u,v)\in\Theta_t$ and $k\in\mathbf Z$, consider a matrix in $\rit{SL}_2(\mathbf Z)$ 
\[
B(t,u,v,k)=\begin{pmatrix}u&(uv-1)/t+uk\\t&v+tk\end{pmatrix}.
\]
 We denote by $\M_{\Theta_t}$ the set of matrices 
$$\{B(t,u,v,k)~|~(u,v)\in\Theta_t,~k \rit{ mod } N/(t^2,N)\}.$$
\begin{lemma}
\begin{enumerate}
\itemx i~The set of matrices $\displaystyle \underset{t|N}{\cup}\M_{\Theta_t}$ is a transversal of the coset decomposition of $\rit{SL}_2(\mathbf Z)$ by $\Gamma_0(N)$.\\
\itemx {ii}~The set of matrices $\{M_\lambda B~|~\lambda\in\SS_N,B\in\underset{t|N}{\cup}\M_{\Theta_t}\}$ is a transversal of the coset decomposition of $\rit{SL}_2(\mathbf Z)$ by $\Gamma_1(N)\{\pm E_2\}$.
\end{enumerate}
\end{lemma}
\begin{proof} The number of elements of the set is $\displaystyle \sum_{t|N}\frac N{(t^2,N)}\varphi((t,N/t))$. This is equal to $[\rit{SL}_2(\mathbf Z):\Gamma_0(N)]$ (see Exercises 11.9 \cite{C1}). It is easy to see that any distinct matrices in the set $\displaystyle \underset{t|N}{\cup}\M_{\Theta_t}$ are not in the same coset. Thus we have (i). The assertion (ii) is obvious from (i).
\end{proof}
Let $\ell_t$ be an integer prime to $t$ and $\ell_t^*$ an integer such that $\ell_t\ell_t^*\equiv 1\mod t$. For the set $\Theta_t$, put 
\[
\ell_t\Theta_t=\{(\ell_t^*u,\ell_t v)|(u,v)\in\Theta_t\}.
\]
Then obviously the set of matrices $\underset{t|N}{\cup}\M_{\ell_t\Theta_t}$ is also a transversal of the coset decomposition.
For an integer $s$ not congruent to $0 \mod N$, let 
$$\phi_s(\tau)=\frac 1{(2\pi i)^2}\wp \left(\frac s N;L_\tau\right)-1/12.
$$
Put $q=\exp (2\pi i\tau/N)$ and $\zeta=\exp(2\pi i/N)$. To consider the $q$-expansion of the function $\phi_s[B(t,u,v,k)]_2$, for an integer $s$, we define two integers $\{s\}$ and $\mu (s)$ by the following conditions:
\[
\begin{split}
&0\le \{s\}\le \frac N2,\quad \mu (s)=\pm 1,\\
&\begin{cases}\mu(s)=1\qquad &\text{if } s\modx {0,N/2}N,\\
             s\equiv \mu (s)\{s\} \mod N\qquad&\text{otherwise.}
\end{cases}
\end{split}
\]
By Lemma 1 of \cite{II}, we have, with $s^*=\mu (st)s(v+tk)$,
{\small
\begin{equation}
\phi_s[B(t,u,v,k)]_2=\phantom{qq\hspace{9.5cm}q}
\end{equation}
}
{\small
\begin{equation*}
\begin{cases}\displaystyle
\frac{\zeta^{s^*}}{(1-\zeta^{s^*})^2}-\sum_{m=1}^{\infty}\sum_{n=1}^{\infty}n(1-\zeta^{s^*n})(1-\zeta^{-s^*n})q^{mnN}&\text{if }\{st\}=0,\\
\phantom{qqqqqqqqqqq}\\
\displaystyle\sum_{n=1}^{\infty}n\zeta^{s^*n}q^{\{st\}n}\\
\displaystyle\phantom{q}-\sum_{m=1}^{\infty}\sum_{n=1}^{\infty}n(1-\zeta^{s^*n}q^{\{st\}n})(1-\zeta^{-s^*n}q^{-\{st\}n})q^{mnN}&\text{otherwise.}\end{cases}
\end{equation*}
}
In particular we note the function $\phi_s[B(t,u,v,k)]_2\in\mathbf Q(\zeta)[[q]]$. 

For an integer $\ell$ prime to $N$, let $\sigma_\ell$ be the automorphism of $\mathbf Q(\zeta)$ over $\mathbf Q$ defined by $\zeta^{\sigma_\ell}=\zeta^\ell$. On a function $f=\displaystyle\sum_m a_mq^m$ with $a_m\in\mathbf Q(\zeta)$, $\sigma_\ell$ acts by $f^{\sigma_\ell}=\displaystyle\sum_m a^{\sigma_\ell}_mq^m$.
\begin{lemma}\label{lem1} Let $\ell$ be an integer prime to $N$ and $\ell^*$ an integer such that $\ell\ell^*\equiv 1\mod N$. Then for $(u,v)\in\Theta_t$ and $k\in\mathbf Z$,
\[
\phi_s[B(t,u,v,k)]_2^{\sigma_\ell}=\begin{cases}\phi_{ls}[B(t,u,v,k)]_2~~&\text{if }\{st\}=0,\\
\phi_s[B(t,\ell^*u,\ell v,\ell k)]_2~~&\text{if }\{st\}\ne 0
\end {cases}
\]
\end{lemma}
\begin{proof}
The $q$-expansion of $\phi_s[B(t,u,v,k)]_2^{\sigma_\ell}$ is given by substituting $s^*$ by $\ell s^*$ in the equation (1). If $\{st\}=0$, then we see $\ell s^*=(\ell s)^*$. If $\{st\}\ne 0$, then $\ell s^*=\mu(st)\ell s(v+tk)=\mu(st)s(\ell v+\ell tk)$. By comparing the $q$-expansion of $\phi_{ls}[B(t,u,v,k)]_2$ or $\phi_s[B(t,\ell^*u,\ell v,\ell k)]_2$ in each case, we have our assertion.
\end{proof}
We consider two subsets $\E_1$ and $\E_2$ of $\E$ given by
\[
\begin{split}
\E_1&=\{\a\in\E~|~(a_1a_2a_3,N)=1\},\\
\E_2&=\{\a\in\E_1~|~(a_i\pm a_3,N)=1 \text{ for }i=1,2\}.
\end{split}
\]
It is noted $\E_1\ne\emptyset$ for $N\geq 7$ (resp.$10$) if $N$ is odd (resp.even) and $\E_2\ne\emptyset$ for $N$ such that $(N,6)=1,N\ge 7$. Further if $N$ is a prime number and $N\ge 7$, then $\E_1=\E_2=\E$.
\begin{exam} Let $\a_1=[2,3,1],\a_2=[2,5,1],\a_3=[1,(N-3)/2,(N-1)/2]$. If $N$ is a positive integer such that $(N,6)=1,N\geq 7$. Then $\a_1,\a_3 \in\E_2$. Further  if $(N,30)=1$, then $\a_2\in \E_2$. The functions $T_{\a_i}$ and $T_{\a_1,\a_2}$ are not constant. See Proposition \ref{progen1}.
\end{exam}
\begin{prop}\label{pro2} Let $\ell$ be an integer prime to $N$ and $\ell^*$ an integer such that $\ell\ell^*\equiv 1\mod N$. Further let $(u,v)\in\Theta_t$ and $k\in\mathbf Z$.
\begin{enumerate}
\itemx i For $\a=[a_1,a_2,a_3]\in\E_1$, we have 
\[
(W_\a\circ B(t,u,v,k))^{\sigma_\ell}=\begin{cases}W_{\ell\a}\circ B(t,u,v,k)~~&\text{if }t=N,\\
W_\a\circ B(t,\ell^* u,\ell v,\ell k)~~&\text{if }t\ne N,
\end{cases}
\]
where $\ell\a=[\{\ell a_1\},\{\ell a_2\},\{\ell a_3\}]$.
\itemx {ii} For a tuple $\A=[\mathfrak a_1,\dots,\mathfrak a_n]~(\a_i\in\E_1)$, we have
\[
(T_{\A,F}\circ B(t,u,v,k))^{\sigma_\ell}=\begin{cases}T_{\A,F}\circ B(t,u,v,k)~~&\text{if }t=N,\\
T_{\A,F}\circ B(t,\ell^* u,\ell v,\ell k)~~&\text{if }t\ne N.
\end{cases}
\]
\end{enumerate}
\end{prop}
\begin{proof} By definition of $W_\a$ we have
\[
W_\a(\tau)=\frac{\phi_{a_1}(\tau)-\phi_{a_3}(\tau)}{\phi_{a_2}(\tau)-
\phi_{a_3}(\tau)}.
\]
Therefore, (i) follows from Lemma~\ref{lem1} and (ii) is obvious from (i) and Proposition~\ref{pro1}.
\end{proof}
It is noted that for $t=1,N$ to obtain the results in Proposition \ref{pro2}, we do not need the condition $\a_i\in\E_1$.
\begin{prop}\label{rationality}
For $\A$ with $\a_i\in\E$, $T_{\A,F}$ and $T_{\A,F}\circ B(1,1,1,-1)$ have Fourier coefficients in $\mathbf Q$.
\end{prop}
\begin{proof}
Since $B(N,u,v,k)\in\Gamma_0(N)$, by Proposition \ref{pro2} (ii), $T_{\A,F}^{\sigma_\ell}=T_{\A,F}$. By the same proposition, we have $(T_{\A,F}\circ B(1,1,1,-1))^{\sigma_\ell}= T_{\A,F}\circ B(1,\ell^*,\ell,-\ell)$. Since $B(1,1,1,-1)B(1,\ell^*,\ell,-\ell)^{-1}\in\Gamma_0(N)$, we see $(T_{\A,F}\circ B(1,1,1,-1))^{\sigma_\ell}=T_{\A,F}\circ B(1,1,1,-1)$.
\end{proof}

Put
\[
\begin{split}
\Phi[W_\a](X,j)=X^{\Psi_1(N)}+\sum_{i=1}^{\Psi_1(N)}C[\a]_iX^{\Psi_1(N)-i},\\
\Phi_{\A,F}(X,j)=X^{\Psi_0(N)}+\sum_{i=1}^{\Psi_0(N)}C_{\A,i}X^{\Psi_0(N)-i},
\end{split}
\]
where $\Psi_0(N)=[\rit{SL}_2(\mathbf Z):\Gamma_0(N)]=\displaystyle N\prod_{p|N}\left(1+\frac 1p\right),~\Psi_1(N)=[\rit{SL}_2(\mathbf Z):\Gamma_1(N)]=\frac{\varphi(N)\Psi_0(N)}2$ and $p$ are prime divisors of $N$.
\begin{thm}\label{th1} 
\begin{enumerate}
\itemx i If $\a\in\E_1$, then the modular equation $\Phi[W_\a]\in\mathbf Q[j][X]$. Further if $N$ is odd and $\a\in\E_2$, then $\Phi[W_\a]\in\mathbf Z[j][X]$.
\itemx {ii}~Let $\A=[\mathfrak a_1,\dots,\mathfrak a_n]$. If $\a_k\in\E_1$ for all $k$, then the modular equation $\Phi_{\A,F}\in\mathbf Q[j][X]$. Further if $N$ is odd and $\a_k\in\E_2$ for all $k$, then $\Phi_{\A,F}\in\mathbf Z[j][X]$.
\end{enumerate}
\end{thm}
\begin{proof}
We know the coefficients $C[\a]_i,~C_{\A,i}\in\mathbf Q(\zeta)((q))$. To show (i), we have only to prove that they are invariant under the action of $\sigma_\ell$ for all $\ell$ prime to $N$. By (i) of Proposition~\ref{pro1}, we see $W_\a\circ (M_\lambda B)=W_{\lambda\a}\circ B$. Thus by Proposition~\ref{pro2}, we have
\[
(W_\a\circ (M_\lambda B(t,u,v,k)))^{\sigma_\ell}=\begin{cases}W_{\a}\circ (M_{\overline{\ell}\lambda}B(t,u,v,k))~&\text{if }t=N,\\
W_\a\circ(M_\lambda  B(t,\ell^* u,\ell v,\ell k))~&\text{if }t\ne N,
\end{cases}
\]
where $\overline{\ell}$ is the element of $\SS_N$ induced by $\ell$. Since $C[\a]_i$ is an elementary symmetric polynomial of $W_\a\circ (M_\lambda B(t,u,v,k))$, we know that $C[\a]_i^{\sigma_\ell}=C[\a]_i$. Therefore We have $C[\a]_i\in\mathbb Q[j]$. Assume that $N$ is odd. Let us consider the $q$-expansions of the function $\phi_a[B]_2-\phi_b[B]_2$ for $a,b\in\mathbf Z,(ab(a\pm b),N)=1$ and $B\in\M_{\Theta_t}$. First of all, let $t\ne N$. Then $
\{at\}\ne\{bt\}$. Let $l=\min(\{at\},\{bt\})$. Then by (1), for an integer $s$
\[
 \phi_a[B]_2(\tau)-\phi_b[B]_2=\pm\zeta^{s}q^l+O(q^{l+1})\in\mathbf Z[\zeta][[q]].
\]
Thus, $W_\mathfrak a\circ B\in\mathbf Z[\zeta]((q))$.  Next we shall consider the case $t=N$. We can take $M_{\Theta_N}=\left\{\bigl(\begin{smallmatrix}
	  1&0\\ N&1
	\end{smallmatrix} \bigr)\right\}$. Put $B=\bigl(\begin{smallmatrix}
	  1&0\\ N&1
	\end{smallmatrix}\bigr)$. By (1), we see
\[
\begin{split}
\phi_a[B]_2(\tau)&-\phi_b[B]_2\\
                &=\frac{\zeta^a(1-\zeta^{b-a})(1-\zeta^{b+a})}{(1-\zeta^a)^2(1-\zeta^b)^2}\\
                &\quad-\sum_{m=1}^\infty \sum_{n=1}^\infty n\{(1-\zeta^{an})(1-\zeta^{-an})-(1-\zeta^{bn})(1-\zeta^{-bn})\}q^{mnN}.
                \end{split}
\]
Let        
\[
\begin{split}
\theta_{a,b}&=\frac{\zeta^a(1-\zeta^{b-a})(1-\zeta^{b+a})}{(1-\zeta^a)^2(1-\zeta^b)^2},\\
h(q)&=-\sum_{m=1}^\infty \sum_{n=1}^\infty n\{(1-\zeta^{an})(1-\zeta^{-an})-(1-\zeta^{bn})(1-\zeta^{-bn})\}q^{mnN}.
\end{split}
\]
Then
\[
\phi_a[B]_2-\phi_b[B]_2=\theta_{a,b}\big(1-\frac{1}{\theta_{a,b}}h(q)\big).
\]
Since $\displaystyle\frac{1-\zeta^s}{1-\zeta^r}\in\mathbf Z[\zeta]^\times$ for integers $r,s$ such that $(rs,N)=1$, we see
        $$\frac{1}{\theta_{a,b}}=\frac{(1-\zeta^a)(1-\zeta^b)}{\zeta^a(1-\zeta^{b-a})(1-\zeta^{b+a})}(1-\zeta^a)(1-\zeta^b)\in\mathbf Z[\zeta].$$
Therefore for some $h(q),f(q)\in\mathbf Z[\zeta][[q]]$
\[
W_\mathfrak a\circ B=\frac{\theta_{a_1,a_3}(1-h(q))}{\theta_{a_2,a_3}(1-f(q))}=\frac{\theta_{a_1,_3}}{\theta_{a_2,_3}}(1+f(q)+f(q)^2+\cdots)(1-h(q)).
\]
Since        
\[
\frac{\theta_{a_1,_3}}{\theta_{a_2,_3}}=\frac{\zeta^{a_1}}{\zeta^{a_2}}\left(\frac{1-\zeta^{a_2}}{1-\zeta^{a_1}}\right)^2\frac{(1-\zeta^{a_3-a_1})(1-\zeta^{a_3+a_1})}{(1-\zeta^{a_3-a_2})(1-\zeta^{a_3+a_2})}\in\mathbf Z[\zeta]^\times,
\]
$W_\mathfrak a\circ B\in\mathbf Z[\zeta][[q]]$. Therefore by (i) of Proposition~\ref{pro1}, we have $W_\mathfrak a\circ (M_\lambda B)\in\mathbf Z[\zeta]((q))$ for all $\lambda\in\SS_n$ and $\displaystyle B\in\cup_{t|N}\Theta_t$. Thus $C[\a]_i\in\mathbf Z[\zeta]((q))$. By applying the above argument, we have $C[\a]_i\in\mathbf Z[j]$. This shows (i). Next we shall prove (ii). By (ii) of Proposition~\ref{pro2}, we have
$$\{(T_{\A,F}\circ B)^{\sigma_\ell}~|~ B\in\M_{\Theta_t}\}=\{T_{\A,F}\circ B~|~ B\in\M_{\ell\Theta_t}\}.
$$
Since $\underset{t|N}{\cup}\M_{\ell\Theta_t}$ is a transversal of coset decomposition of $\rit{SL}_2(\mathbf Z)$ by $\Gamma_0(N)$, we obtain $C_{\A,i}^{\sigma_\ell}=C_{\A,i}$. This shows $C_{\A,i}\in\mathbf Q[j]$. If $N$ is odd and $\a\in\E_2,~\lambda\in\SS_N$, then $\lambda\a\in\E_2$. Proposition~\ref{pro1} shows $T_{\A,F}\circ B\in\mathbf Z[\zeta]((q))$. Therefore $C_{\A,i}\in\mathbf Z[\zeta]((q))$. Since $C_{\A,i}^{\sigma_\ell}=C_{\A,i}$, this shows $C_{\A,i}\in\mathbf Z[j]$.
\end{proof}
Let $K$ be an imaginary quadratic field and $\K_N$ the ray class field modulo $N$ over $K$.  
\begin{thm}\label{singv} Assume that $N$ is odd. Let $\alpha$ be an element of $\H$ such that $K=\mathbf Q(\alpha)$. 
\begin{enumerate}
\itemx i If $\a\in\E_2$, then $W_\a(\alpha)$ is a unit of $\K_N$.\\
\itemx {ii} ~Let $\A=[\mathfrak a_1,\dots,\mathfrak a_n]$. If $\a_k\in\E_2$ for all $k$, then $T_{\A,F}(\alpha)$ is an algebraic integer of $\K_N$.
\end{enumerate}
\end{thm}
\begin{proof}
By Complex multiplication theory, $j(\alpha)$ is an algebraic integer. Theorem~\ref{th1} shows that $\Phi[W_\a](X,j(\alpha))$ and $\Phi_{\A,F}(X,j(\alpha))$ are monic polynomials with algebraic integer coefficients. Thus $W_\a(\alpha),T_{\A,F}(\alpha)$ are algebraic integers. By Corollary to Theorem 2 in \S 10.1 of \cite{LA}, they are in $\K_N$. Let $\a'=[a_2,a_1,a_3]$. Since $W_\a^{-1}=W_{\a'}$ and $\a'\in\E_2$, $W_\a(\alpha)^{-1}$ is an algebraic integer. Hence it is a unit. 
\end{proof}
\section{Ray class field and ring class field}
Let $K$ be a subfield of $\mathbf C$ and $\Gamma$  a subgroup of $\rit{SL}_2(\mathbf Z)$ of finite index. We denote by $A(\Gamma)_K$ the field of all modular functions with respect to $\Gamma$ having Fourier coefficients in $K$. Further put $A_0(N)_K=A(\Gamma_0(N))_K,~A_1(N)_K=A(\Gamma_1(N))_K$. Let $\zeta=\exp(2\pi i/N)$.
\begin{prop}\label{progen}
Put $\a_1=[2,3,1],\a_2= [2,5,1]$. If $N\ge 11,N\ne 12$, then 
\[
A_1(N)_{\mathbf Q(\zeta)}=\mathbf Q(\zeta)(j,W_{\a_1})=\mathbf Q(\zeta)(j,W_{\a_2})=\mathbf Q(\zeta)(W_{\a_1},W_{\a_2}).\]
\end{prop}
\begin{proof}
The assertion is deduced from the result $A_1(N)_{\mathbf C}$ $=\mathbf C(j,W_{\a_i})$ $=\mathbf C(W_{\a_1},W_{\a_2})$ and $W_{\a_i}\in A_1(N)_{\mathbf Q(\zeta)}(i=1,2)$ in Lemma 1 and Theorems 1 and 5 of \cite{II}.
\end{proof}
Let $m$ and $n$ be non-negative integers.  Put $F=X_1^mX_2^n$ and $\A=[\a_1,\a_2]$ with $\a_1=[2,3,1],\a_2= [2,5,1]$. For a while we shall consider the function $T_{\A,F}$. By Theorem 3.2 of \cite{I1}, for any $\b=[b_1,b_2,b_3]\in\E$, the order of the $q$-expansion of $W_\b$ at the point $u/t$ is equal to $\min(\{tb_1\},\{tb_3\})-\min(\{tb_2\},\{tb_3\})$. In particular, the order of $q$-expansion of $W_\b\circ B(t,u,v,k)$ depends only on $t$ and it equals to that of $W_\b$ at the point $1/t$. For any integers $a,b$ and $c$, we see $\{\{ab\}c\}=\{a\{bc\}\}$. Thus the order of $q$-expansion of $W_{\lambda\a_i}\circ B(t,u,v,k)$ is that of $W_{\a_i}$ at the point $1/\{\lambda t\}$. Let $\omega_i(\ell)$ be the order of $q$-expansion of $W_{\a_i}$ at the point $1/\ell$ for $\ell\in\mathbf Z,~1\leq \ell\leq N/2$. By \S 3 of \cite{II}, we know $\omega_i(\ell)<0$ if and only if $\ell>\frac{2N}5$ (resp. $\frac{3N}7$) for $i=1$ (resp. $i=2$). We have $\omega_{i}(\ell)=(i+1)N-(2i+3)\ell$ for $2N/5<\ell\leq N/2$ and ,in this range, obviously $\omega_{i}(\ell)$ is a strictly decreasing function of $\ell$. Furthermore $\omega_{i}(\ell)\modx 0{(\ell,N)}$. 
\begin{prop}\label{progen1} Assume that $N$ is a prime number and $N>7$. Put $\a_1=[2,3,1],\a_2= [2,5,1]$ and $\a_3=[1,(N-3)/2,(N-1)/2]$. Then for $i=1,3$ and $j=1,2,3$
\[A_0(N)_{\mathbf Q}=\mathbf Q(T_{\a_i},T_{\a_2})=\mathbf Q(T_{\a_j},T_{\a_1,\a_2}).\] 
\end{prop}
\begin{proof}
  Put $T_i=T_{\a_i}$ for $i=1,2,3$ and $T_4=T_{\a_1,\a_2}$. 
Since $N$ is a prime number, the group $\Gamma_0(N)$ has two cusps represented by $i\infty$ and $1$. By Theorem 3.2 of \cite{I1}, for any $\b=[b_1,b_2,b_3]\in\E$, $W_\b$ is regular at the point $i\infty$. Therefore the functions $T_i$ $(i=1,\dots,4)$ are regular at $i\infty$. Let us denote by $d_i$  the order of the pole of $T_i$ at the  cusp $1$. We know $\omega_{i}(\lambda)$ has the smallest value only for $\lambda=(N-1)/2$. Thus, we have $d_1=(N-5)/2,d_2=(N-7)/2$ and $d_4=N-6$. Let us determine $d_3$. The function $W_{\a_3}$ has a pole of order $(N-5)/2$ at $1$. Let $\lambda >1$. The function $W_{\lambda\a_3}$ has a pole  at $1$ if $\lambda<\{\lambda(N-1)/2\}<\{\lambda(N-3)/2\}$ or $\lambda<\{\lambda(N-3)/2\}<\{\lambda(N-1)/2\}$. In the former case, the order $d_\lambda$ of pole of $W_{\a_3}$ at $1/\lambda$ is $\{\lambda(N-1)/2\}-\lambda$. Since $\{\lambda(N-1)/2\}<\{\lambda(N-3)/2\}$, we know $\{\lambda(N-1)/2\}\leq(N-3)/2$. Thus $d_\lambda<(N-5)/2$. In the latter case, $d_\lambda=\{\lambda(N-3)/2\}-\lambda$. Since $\lambda>1, \{\lambda(N-3)/2\}\leq (N-3)/2$, we know $d_\lambda<(N-5)/2$. Therefore we have $d_3=(N-5)/2$.  Proposition \ref{rationality} shows that $T_i\in A_0(N)_{\mathbf Q}$. Since the modular curve $X_0(N)$ of $\Gamma_0(N)$ is defined over $\mathbf Q$, by Proposition 2.6 (a) in Chapter II of \cite{SJ}, $d_i=[A_0(N)_{\mathbf Q}:\mathbf Q(T_i)]$. Since $((N-5)/2,(N-7)/2)=1$ and $((N-5)(N-7),(N-6))=1$, we have our assertion.\end{proof}
\begin{thm}\label{progen3}Let $m$ and $n$ be non-negative integers. Assume that $N$ does not divide $5m+7n$ {\rm (resp. $2(5m+7n)$)} and $N>9$ {\rm (resp. $36$)} in the case $N$ is odd {\rm(resp. even)}. Put $\A=[\a_1,\a_2]$ and $F=X_1^mX_2^n$. Further assume that $N\not\equiv 0\mod 4$ in the case $m+n$ is even. Then we have $A_0(N)_{\mathbf Q}=\mathbf Q(j,T_{\A,F})$.
\end{thm}
\begin{proof}
Put $T=T_{\A,F}$. By Theorem 3 of Chapter  6 of \cite{LA}, the field $A(\Gamma (N))_{\mathbf Q(\zeta)}$ is a Galois extension over $\mathbf Q(j)$ with the Galois group $\rit{GL}_2(\mathbf Z/N\mathbf Z)/\{\pm E_2\}$ and the field $A_0(N)_{\mathbf Q}$ is the fixed field of the subgroup $\displaystyle\left\{\begin{pmatrix}*&*\\0&*\end{pmatrix}\right\}/\{\pm E_2\}$. Since $T\in A_0(N)_{\mathbf Q}$, to prove the assertions, it is sufficient to show that if $T\circ A=T$ for $A\in\rit{SL}_2(\mathbf Z)$, then $A\in\Gamma_0(N)$. Let us consider the transversal $\{B(t,u,v,k)\}$ of the coset decomposition of $\rit{SL}_2(\mathbf Z)$ by $\Gamma_0(N)$. Let $\omega(\ell)$ be the order of $q$-expansion of $W_{\a_1}^mW_{\a_2}^n$ at the point $1/\ell$. Then obviously $\omega(\ell)=m\omega_1(\ell)+n\omega_2(\ell)$. Let $t$ be a divisor of $N$. If $\lambda$ runs over $\SS_N$, then $\{\lambda t\}$ runs over all integers $u$ such that $0\leq u\leq N/2,~(u,N)=t$. Therefore $d\geq\min\{\omega(\ell)~|~0\leq\ell\leq N/2,~(\ell,N)=u\}$. Furthermore if $\omega(\ell)$ has the smallest value for only one $\ell$, then we have equality. Let $u_t$ be the greatest integer such that $(u_t,N)=t$ and $u_t\leq N/2$. 
 Let $t\ne N$.  Assume that  $T\circ B(t,u,v,k)=T$. Put $L=\displaystyle B(1,1,1,-1)=\begin{pmatrix}1&-1\\1&0\end{pmatrix}$. Then $T\circ (B(t,u,v,k) L)=T\circ L$. We know $\displaystyle B(t,u,v,k)L=\begin{pmatrix}*&*\\t(k+1)+v&-t\end{pmatrix}$.  Let $\delta=(t(k+1)+v,N)$. Then we can take an integer $\xi$ so that $\xi ((k+1)t+v)+\delta t\modx 0N$ and $(\xi,\delta)=1$. For an integer $\eta$ such that $\xi\eta\modx 1\delta$, put $\displaystyle A=\begin{pmatrix}\eta&(\xi\eta-1)/\delta\\ \delta&\xi\end{pmatrix}$. Since $B(t,u,v,k) LA^{-1}\in\Gamma_0(N)$, we have $T\circ A=T\circ L$. Let $d$ be the order of $q$-expansion of $T\circ A$ and $d_1$ the order of $q$-expansion of $T\circ L$. In particular, the assumption implies that $d=d_1\modx 0\delta$. In the case $\delta\ne 1$, we shall show that $d\ne d_1$. If $N$ is even, then $u_{N/2}=N/2$. If $\delta\ne N/2$, then $u_\delta$ is as follows.
\begin{table}[h]
\begin{center}
\begin{tabular}{|c|c|c|c|}
\hline
$N~ \rit{mod} 4$&$N/\delta~ \rit{mod} 4$&$u_1$&$u_\delta$\\\hline
$1,3$&$1,3$&$(N-1)/2$&$(N-\delta)/2$\\\hline
$2$&$0,2$&$(N-4)/2$&$(N-4\delta)/2$\\\hline
$2$&$1,3$&$(N-4)/2$&$(N-\delta)/2$\\\hline
$0$&$1,3$&$(N-2)/2$&$(N-\delta)/2$\\\hline
$0$&$2$&$(N-2)/2$&$(N-4\delta)/2$\\\hline
$0$&$0$&$(N-2)/2$&$(N-2\delta)/2$\\\hline
\end{tabular}
\end{center}
\end{table}
\newpage
If we put $u_1=(N-\epsilon)/2$ with $\epsilon=1$ (resp. $2,4$) in the case $N$ is odd (resp. even), we see easily $d_1=\omega(u_1)=((5m+7n)\epsilon-(m+n)N)/2$ and $d\geq\min(0,\omega(u_\delta))$. It is noted that our assumption implies $d_1<0$. If $\delta=N$, then $d\geq 0$. Thus $d\ne d_1$. If $\delta=N/2$, then $d\modx 0{N/2}$. By assumption, $d_1\not\equiv 0\mod N/2$. This implies $d\ne d_1$. Let $\delta\ne 1,N/2,N$.
Except the case $N\modx 24$ and $\delta=2$, we have $u_\delta<u_1$. Thus $d\ne d_1$. In the exceptional case, we have $u_\delta>u_1$. Since there exists only one $\lambda$ such that $\{2\lambda\}=u_{2}$, we have $d<d_1$. Let us consider the case $\delta=1$. Then $d=d_1$. By (1), for a matrix $\displaystyle M=\begin{pmatrix}*&*\\1&k\end{pmatrix}$ of $\rit{SL}_2(\mathbf Z)$ and $0<s\leq N/2$, we have
\[
\varphi_s\circ M=\zeta^{s^*}q^{s}+\zeta^{-s^*}q^{N-s}+2\zeta^{2s^*}q^{2s}+2\zeta^{-2s^*}q^{2(N-s)}-q^N+\text{(higher terms)},
\]
where $s^*=\mu(s)sk=sk$. If we put $s=s_r=\{ru_1\}$ for $r=1,2,3,5$, then 
$s_r=(N-r\epsilon)/2,s_r^*=ru_1k$ for $r=1,3,5$ and $s_2=\epsilon,s_2^*=-2u_1k$. Since $(N-\epsilon)/2>2\epsilon$, we have
\begin{equation}\label{eq2}
\begin{split}
(\varphi_{s_2}-\varphi_{s_1})\circ M&=\zeta^{-2u_1k}q^\epsilon(1+2\zeta^{-2u_1k}q^\epsilon +O(q^{\epsilon+1})),\\
(\varphi_{s_3}-\varphi_{s_1})\circ M&=\zeta^{3u_1k}q^{(N-3\epsilon)/2}(1-\zeta^{-2u_1k}q^\epsilon+O(q^{\epsilon+1})),\\
(\varphi_{s_5}-\varphi_{s_1})\circ M&=\zeta^{5u_1k}q^{(N-5\epsilon)/2}(1+O(q^{\epsilon+1})),
\end{split}
\end{equation}
where the notation $O(q^n)$ denotes a $q$-series of order greater than or equal to $n$. Because the assumption for $N$ implies $d_1+\epsilon<0, \omega(u_1-1)$, we see by \eqref{eq2}, 
\[
T\circ M=\zeta^{-(5m+7n)u_1k}q^{d_1}(1+(3m+2n)\zeta^{-2u_1k}q^\epsilon+O(q^{\epsilon+1})).\]
If we compare the coefficients of $T\circ A~(k=\xi)$ with those of $T\circ L (k=0)$, we see $\zeta^{-(5m+7n)u_1\xi}=\zeta^{-2u_1\xi}=1$. If $(5m+7n)$ is odd, then, since $(u_1,N)=1$, we have $\xi\modx 0N$. Since $\xi(t(k+1)+v)+t\modx 0N$, we have $t\modx 0N$. This gives a contradiction. Obviously if $(5m+7n,N)=1$, we have also a contradiction.  For the case $5m+7n$ is even, we have $\zeta^{-2u_1\xi}=1$. This shows $2t\modx 0N$. Therefore if $N$ is odd, we have a contradiction. Let $N\modx 24$ and $t=N/2$. It is noted $\Theta_{N/2}=\{B(N/2,1,1,k)~|~k=0,1\}$. Since $(N/2,2)=1$, we can take integers $x$ and $y$ such that $(N/2)x+2y=1$. Consider a matrix $\displaystyle A=\begin{pmatrix}x&-1\\2y&t\end{pmatrix}$ of $\rit{SL}_2(\mathbf Z)$. 
It is easy to see that $B(N/2,1,1,k)AB(1,1,1,k(N/2)^2-1)^{-1}$, $AB(2,1,1,-y)^{-1}\in\Gamma_0(N)$. Therefore we have $T\circ B(1,1,1,k(N/2)^2-1)=T\circ B(2,1,1,-y)$. However the above argument for $N\modx 24$ and $\delta=2$ shows the order of $q$-expansions of the functions $T\circ B(1,1,1,k(N/2)^2-1)$ and $T\circ B(2,1,1,-y)$ are distinct. 
\end{proof}
\begin{thm}\label{generator}
Let $\alpha\in\H$ such that $\mathbf Z[\alpha]$ is a maximal order of $K$. Further let $\a_1=[2,3,1]$ and $\a_2=[2,5,1]$. If $N=11$ or $N\geq 13$, then 
\[
\K_N=K(\zeta,j(\alpha),W_{\a_1}(\alpha))=K(\zeta,j(\alpha),W_{\a_2}(\alpha))=K(\zeta,W_{\a_1}(\alpha),W_{\a_2}(\alpha)).\]
\end{thm}
\begin{proof}
Our assertion follows from Theorems 1 and 2 of \cite{GA} and Proposition~\ref{progen}.
\end{proof}
For a positive integer $m$, let $O_m$ be the order of conductor $m$ of $K$ and $R_m$ the ring class field associated with the order $O_m$. Consider the group 
\[
\Gamma^0(N)=\left\{\left.\begin{pmatrix}a&b\\c&d\end{pmatrix}\in SL(2,\mathbf Z)\right|b\equiv 0 \mod N\right\}.
\]
\begin{prop}\label{provalue}
Let $\theta\in\H$ such that $O_f=\mathbf Z[\theta]$ and $f_\theta(X)=X^2+BX+C~(B,C\in\mathbf Z)$ the minimal polynomial of $\theta$. 
\begin{enumerate}
\itemx i If $h\in A_0(N)_{\mathbf Q}$ and $h$ is pole-free at $\theta$, then $h(\theta)\in R_{fN}$.\\
\itemx {ii} If $h\in A(\Gamma^0(N))_{\mathbf Q}$, $h$ is pole-free at $\theta$ and $N|C$,then $h(\theta)\in R_{f}$.
\end{enumerate}
\end{prop}
\begin{proof}
Let us use the notation in \S 2 of \cite{GA}. For a prime number $p$, consider groups 
\[
\begin{split}
 U_p=\left\{\left.\begin{pmatrix}a&b\\c&d\end{pmatrix}\in\rit{GL}_2(\mathbf Z_p)~\right |c\in N\mathbf Z_p\right\},\\
 V_p=\left\{\left.\begin{pmatrix}a&b\\c&d\end{pmatrix}\in\rit{GL}_2(\mathbf Z_p)~\right |b\in N\mathbf Z_p\right\}.
\end{split}
\]
Put $\displaystyle U=\prod_p U_p,~V=\prod_p V_p$. Then $U$ (resp.$V$) is the subgroup of $\displaystyle\prod_p\rit{GL}_2(\mathbf Z_p)$ with fixed field $F_0=A_0(N)_{\mathbf Q}$ (resp.$F^0=A(\Gamma^0(N))_{\mathbf Q}$). Let $O=O_f$ be the order of $K$ of conductor $f$. By Theorems 5.5 and 5.7 of \cite{SG}, we have an exact sequence
\[1\longrightarrow O^*\longrightarrow \prod_pO_p^*\longrightarrow \rit{Gal}(K^{ab}/K(j(\theta))\longrightarrow 1.\]
Let $\displaystyle g_\theta=\prod_p(g_{\theta})_p:\prod_pO_p^*\longrightarrow \prod_p\rit{GL}_2({\mathbf Z}_p)$ be the map defined by (4) and (5) in \cite{GA}. Since by the definition, for $s,t\in\mathbf Z_p$, 
\[
(g_{\theta})_p(s\theta+t)=\begin{pmatrix}t-Bs&-Cs\\s&t\end{pmatrix},
\]
we have $\displaystyle (g_\theta)_p^{-1}(U_p)=(\mathbf Z_p^*+N\mathbf Z_p\theta)\cap O_p^*=(O_{fn})_p^*$. Therefore we have $\displaystyle g_\theta^{-1}(U)=\prod_p(O_{Nf})_p^*$. If $N|C$, then $\displaystyle (g_\theta)_p^{-1}(V_p)=O_p^*$ and $\displaystyle g_\theta^{-1}(V)=\prod_pO_p^*$. By class field theory the groups $\displaystyle \prod_p O_p^*$ and $\displaystyle \prod_p (O_{fN})_p^*$ correspond to $\rit{Gal}(K^{ab}/R_{f})$ and $\rit{Gal}(K^{ab}/R_{fN})$ respectively. By Theorem 2 of \cite{GA}, we see $R_{f}=K(F^0(\theta))$ and $R_{fN}=K(F_0(\theta))$. Therefore we have our assertions. 
\end{proof}
\begin{cor}\label{cor1}
Let the notation be the same as in Proposition~\ref{provalue}. Let $\A=[\a_1,\dots,\a_n]$ with $\a_i\in\E_1$ for all $i$. Then we have the followings.
\begin{enumerate}
\itemx i $T_{\A,F}(\theta)\in R_{fN}$.
\itemx {ii}~If $N|C$, then $T_{\A,F}(-1/\theta)\in R_{f}$.
\end{enumerate}
\end{cor}
\begin{proof}
Since $\Gamma^0(N)=S^{-1}\Gamma_0(N)S$ with $S=\big(\begin{smallmatrix}0&-1\\1&0\end{smallmatrix}\big)$, $T_{\A,F}\circ S$ is a modular function with respect to $\Gamma^0(N)$. By the result for $t=N$ in (ii) of Proposition~\ref{pro2}, we know $T_{\A,F}\in\mathbf Q((q))$. Since $B(1,1,1,-1)S^{-1}\in\Gamma_0(N)$, we have $T_{\A,F}\circ S=T_{\A,F}\circ B(1,1,1,-1)$. By Proposition \ref{rationality}, we know  $T_{\A,F}\in A_0(N)_{\mathbf Q}$ and $T_{\A,F}\circ S\in A(\Gamma^0(N))_{\mathbf Q}$. Our assertions follow from Proposition~\ref{provalue}.
\end{proof}
\begin{thm}\label{singt}
Put $\a_1=[2,3,1],\a_2=[2,5,1],\a_3=[1,(N-3)/2,(N-1)/2]$ and $\A=[\a_1,\a_2]$. Put $F=X_1^mX_2^n$ with non-negative integers $m$ and $n$. Let $\theta\in\H$ such that $O_f=\mathbf Z[\theta]$ and $f_\theta(X)=X^2+BX+C~(B,C\in\mathbf Z)$ the minimal polynomial of $\theta$. Then we have followings.
\begin{enumerate}
\itemx {i}~If $N$ is a prime number and $N>7$, then
\[
R_{fN}=K(T_{\a_i}(\theta),T_{\a_2}(\theta))=K(T_{\a_j}(\theta),T_{\a_1,\a_2}(\theta)),\]
for $i=1,3$ and $j=1,2,3$. Further if $N|C$, then
\[
R_{f}=K(T_{\a_i}(-1/\theta),T_{\a_2}(-1/\theta))=K(T_{\a_j}(-1/\theta),T_{\a_1,\a_2}(-1/\theta)),\]
for $i=1,3$ and $j=1,2,3$.
\itemx {ii} Assume that $N$ does not divide $5m+7n$ {\rm (resp. $4(5m+7n)$)} and $N>9$ {\rm (resp. $36$)} in the case $N$ is odd {\rm(resp. even)}. Further assume that $N$ is not divided by $4$ in the case $m+n$ is even. Then $R_{fN}=K(j(\theta),T_{\A,F}(\theta))$. Further if $N|C$, then $R_{f}=K(j(\theta),T_{\A,F}(-1/\theta))$
\end{enumerate}
\end{thm}
\begin{proof}
In the proof of Propositions \ref{provalue} we showed $R_{Nf}=K(F_0(\theta)), R_f=K(F^0(\theta))$. Therefore the assertions follow from Propositions \ref{progen1} and Theorem \ref{progen3}.\end{proof}
\section{Class polynomials of $T_{\A,F}$} 
Let $O$ be the order of conductor $f$ of an imaginary quadratic field $K$. Let $D$ be the discriminant and $C(O)$ the (proper) ideal class group of $O$. We denote by $h(D)$ the class number of $O$. Let $\alpha\in K\cap\H$ and $AX^2+BX+C=0$ be the primitive  minimal equation with integral coefficients of $\alpha$ over $\mathbf Q$. If $D=B^2-4AC$, then we say $\alpha$ is an element of discriminant $D$. We put $I_\alpha=[A,(-B+\sqrt D)/2]=\mathbf Z A+\mathbf Z ((-B+\sqrt D)/2)$. Then $I_\alpha$ is an ideal of $O$. To compute the singular values of the functions $T_{\A,F}$, we use an $N$-system for $O$ introduced by Schertz \cite{RS}.
\begin{defn}
Let $\N$ be a set of $h(D)$ elements $\alpha_i\in K\cap \H$ of discriminant $D$. Let $A_iX^2+B_iX+C_i=0$ be the primitive integral minimal equation of $\alpha_i$ and $I_{\alpha_i}=[A_i,(-B_i+\sqrt D)/2]$. We say $\N$ is an $N$-system for $O$ if following conditions are satisfied:
\begin{enumerate}
 \item[1.] $(A_i,N)=1,~N|C_i, B_i\modx{B_j}{2N}$ for every $i,j$,
 \item[2.] the set of ideals $\{I_{\alpha_i}\}$ is a transversal of $C(O)$.
\end{enumerate}
\end{defn}
Let $\N$ be an $N$-system for $O$. Then by Complex multiplication theory, for each $\alpha_i\in\N$, $j(\alpha_i)$ is an algebraic integer and generates the ring class field $R_f$ associated with the order of conductor $f$ and they are conjugate to each other over $\mathbf Q$ (see \S 11.D of \cite{C1}). For singular values $T_{\A,F}(-1/\alpha_i)$ we have
\begin{thm}\label{th3} Let $N$ be a positive integer such that $\E_2$ is not empty. Put $\A=[\a_1,\dots,\a_n]$ with $\a_i\in\E_2$. Let $\N=\{\alpha_i\}$ be an $N$-system for $O$. Then we have $T_{\A,F}(-1/\alpha_i)\in R_f$ and they are conjugate to each other over $K$.
\end{thm}
\begin{proof}
 Since $\Gamma^0(N)=S^{-1}\Gamma_0(N)S$ with $S=\big(\begin{smallmatrix}0&-1\\1&0\end{smallmatrix}\big)$, $T_{\A,F}\circ S$ is a modular function with respect to $\Gamma^0(N)$. In the proof of Corollary \ref{cor1}, we showed $T_{\A,F}\circ S\in\mathbf Q((q))$. Therefore the assertion follows from Theorem 3.1 of \cite{ES} and Theorem~\ref{singv}.
\end{proof}

For a modular function $g(\tau)$ with respect to $\Gamma_0(N)$ and an $N$-system $\N=\{\alpha_i\}$, we define the class polynomial $H_\N[g](X)$ of $g(\tau)$ by
 
\[
H_\N[g](X)=\prod_{i=1}^{h(D)}\left(X-g(-1/\alpha_i)\right).
\]

The next assertion follows from Theorem~\ref{th3}.
\begin{thm}
Let $O_K$ be the maximal order of $K$. Then the class polynomial $H_\N[T_{\A,F}](X)\in O_K[X]$.
\end{thm}
Let $B$ be an integer such that $B^2\equiv D\pmod {4N}$. Proposition 3 of \cite{RS} shows the existence of $N$-system containing the number $(-B+\sqrt{D})/2$.
By Lemma 3.1 of \cite{YAI}, we know the class polynomials of a modular function $g$ related to $N$-systems depend only on integers $B$, considered mod $2N$. We shall fix an $N$-system containing $(-B+\sqrt{D})/2$ and denote it by $\N_B$. In the followings, we give some examples of modular equations and class polynomials of the functions $f=T_\a$ or $T_{\a_1,\a_2}$. We shall denote by $H_B(X)$ the class polynomial $H_{\N_B}[f]$ in the case the function $f$ is clearly indicated and any confusion can not occur.
{\small
\begin{exam}\label{exam1}
\begin{enumerate}
\itemx 1~Let $N=7,\a=[2,3,1]$. Consider the function $T_\a$. Then the modular equation $\Phi(X,j)$ of $T_\a$ is given by
\[
\begin{split}
\Phi(X&,j)=X^8-36X^7+546X^6-4592X^5+23835X^4\\
&-80304X^3+176050X^2-(j+232500)X+140625+8j.
\end{split}
\]
\begin{enumerate}
\itemx a~Let $D=-3,B=5$. Then $h(-3)=1,\N_{5}=\{(-5+\sqrt{-3})/2\}$. We have the class polynomial $$H_{5}(X)=X-3(1+\sqrt{-3})/2.$$ Thus $T_\a((5+\sqrt{-3})/14)=3(1+\sqrt{-3})/2$. Since $j((1+\sqrt{-3})/2)=0$, we have $\Phi(X,0)=(X^2-3X+9)(X^2-11X+25)^3$. In fact, $T_\a((5+\sqrt{-3})/14)$ is a root of the factor $X^2-3X+9=0$.
\itemx b~ Let $D=-59,B=5$. Then we have $h(-59)=3$ and 
{\footnotesize
\[
\begin{split}
&\N_5=\{(-5+\sqrt{-59})/2, (-5+\sqrt{-59})/6,(23+\sqrt{-59})/6\}\\ 
&H_5(X)=X^3+\frac{15-7\sqrt{-59}}2X^2+\frac{-357+45\sqrt{-59}}2X+\frac{717+\sqrt{-59}}2.
\end{split}
\]
}
\end{enumerate}
\itemx 2~Let $N=13,D=-3,B=7,\a=[5,3,1]$. Take $\N_7=\{(-7+\sqrt{-3})/2\}$. Then the modular equation $\Phi(X,j)$ of $T_\a$ and the value $T_\a(7+\sqrt{-3})/26)$ are given by {\footnotesize
\[ 
\begin{split}                                                       
&\Phi(X,j)=(X^2-9X+27)(X^4-21X^3+167X^2+-604X+848)^3-j(X-7),\\
&T_\a((7+\sqrt{-3})/26)=(9+3\sqrt{-3})/2.
\end{split}
\]
}
Thus in fact $T_{\A,F}((7+\sqrt{-3})/26)$ is a root of $X^2-9X+27=0$D
\itemx 3~Let $N=11,~\a=[2,5,1],D=-7,B=9$. Then $\N_9=\{(-9+\sqrt{-7})/2\}$ and we have $T_\a((9+\sqrt{-7})/44)=(5+\sqrt{-7})/2$ and the modular equation
{\footnotesize
\[\begin{split}
\Phi&(X,j)=X^{12}-84X^{11}+2970X^{10}-57772X^9+680559X^8-5062728X^7\\
&-(22j-24250028)X^6+(561j-75844824)X^5-(2981j-157525071)X^4\\
&-(1177j+217265444)X^3+(26477j+193124250)X^2\\
&-(j^2+31316j+101227452)X+18j^2+4261j+24137569.
\end{split}
\]
}
Since $j((1+\sqrt{-7})/2)=-15^3$, we have 
{\footnotesize
\[
\begin{split}
&\Phi(X,-15^3)=(X^{10}-79X^9+2567X^8-44305X^7+438498X^6-2515798X^5\\
&+8237304X^4-16425295X^3+19561039X^2+15914486X+26848493)\\
&\times(X^2-5X+8).
\end{split}
\]
}
Therefore, we know $T_\a((9+\sqrt{-7})/44)$ is a root of the factor $X^2-5X+8$. \end{enumerate}
\end{exam}
\begin{exam}\label{exam2}
Let $N=11,~\a=[2,3,1],\b=[2,3,5]$.  Consider the function $T_{\a,\b}$. Then we give the coefficients $C_i$ of the modular equation $\displaystyle\Phi(X,j)=X^{12}+\sum_{i=1}^{12}C_iX^{12-i}$ in the table below.\newline
\rit{(1)}~ Let $D=-83,B=7$. Then we have $h(-83)=3$ and
{\footnotesize
\[\begin{split}
&\N_7=\{(-7+\sqrt{-83})/2,(-7+\sqrt{-83})/6,(-29+\sqrt{-83})/6\},\\
&H_7(X)=X^3-(361481+7136\sqrt{-83})X^2+(57020581+25984608\sqrt{-83})X\\
&\phantom{aaaaaaaaaa}+1683573861-404390656\sqrt{-83}.
\end{split}
\]
}
\noindent
\rit{(2)} Let $D=-39,B=7$. Then we have $h(-39)=4$ and 
{\footnotesize
\[
\begin{split}
&\N_7=\{(-7+\sqrt{-39})/2, (-7+\sqrt{-39})/4,(-29+\sqrt{-39})/4,(-51+\sqrt{-39})/8\},\\
&H_7(X)=X^4+(-4720+231\sqrt{-39})X^3+(1491643-329343\sqrt{-39})X^2/2\\
&\phantom{aaaaaa}+(-38934427+9970611\sqrt{-39})X/2+64994911-47480958\sqrt{-39}.
\end{split}
\]
}
\end{exam}
\begin{table}[h]
\begin{tabular}{|l|l|l|l|}
\hline
$i$&{\scriptsize$C_i$}&$i$&\\\hline
{\scriptsize$1$}&{\scriptsize$3660$}&{\scriptsize$2$}&{\scriptsize$4754178$}\\\hline
{\scriptsize$3$}&{\scriptsize$21879j+2517699932$}&{\scriptsize$4$}&{\scriptsize$8917579j+450023862255$}\\\hline
\end{tabular}\newline
\begin{tabular}{|l|l|}\hline
{\scriptsize$5$}&{\scriptsize$10912j^2-21727187108j+28522470464664$}\\\hline
{\scriptsize$6$}&{\scriptsize$18536243j^2+439266301210j+155307879800348$}\\\hline
{\scriptsize$7$}&{\scriptsize$1419j^3+6356028822j^2-4268224633178j-22718073239498472$}\\\hline
{\scriptsize $8$}&{\scriptsize$1663761j^3+70427463557j^2-129554423289764j+430444117263292143$}\\\hline
{\scriptsize$9$}&{\scriptsize$66j^4-100966360j^3+544875974962j^2+1322596244939332j-4047340123195216100$}\\\hline
{\scriptsize$10$}&{\scriptsize$82687j^4-2985616392j^3+3765768493971j^2-9777105305922130j+21981914597781276930$}\\\hline
{\scriptsize$11$}&{\scriptsize$j^5+1956838j^4-26707875453j^3+49826805469384j^2+21725643544520963j$}\\
&{\scriptsize$-67067772106836815988$}\\\hline
{\scriptsize$12$}&{\scriptsize$1229j^5+29053078j^4-41072974661j^3-92728235099098j^2+68572479313531217j$}\\
&{\scriptsize$+93554961663154376449$}\\\hline
\end{tabular}
\end{table}
\begin{exam}\label{exam3}
Let $N=17,\a=[1,2,7],\b=[1,2,3]$. Consider the function $T_{\a,\b}$. Let $D=-84,B=8$. Then we have $h(-84)=4$ and 
\[
\begin{split}
&\N_8=\{-8+\sqrt{-21},(43+\sqrt{-21})/11,(-8+\sqrt{-21})/5,(9+\sqrt{-21})/3\},\\&H_8(X)=X^4+(779-157\sqrt{-21})X^3+(-41194-175\sqrt{-21})X^2\\
&\phantom{aa}+(690208+81256\sqrt{-21})X-3246464-566976\sqrt{-21}.
\end{split}
\]
\end{exam}
}

\vspace{5mm}

{\small
\begin{tabular}{ll}
Faculty of Liberal arts and Sciences \phantom{aaaaaa}& Graduate School of Science\\
Osaka Prefecture University &Osaka Prefecture University \\
1-1 Gakuen-cho, Naka-ku Sakai &1-1 Gakuen-cho, Naka-ku Sakai\\
 Osaka, 599-8531 Japan& Osaka, 599-8531 Japan \\
e-mail:\quad ishii@las.osakafu-u.ac.jp &
\end{tabular}
}

\begin{thebibliography}{99}
\bibitem{C1} D.Cox, Primes of the form $x^2+ny^2$, A Wiley-Interscience Publication, John Wiley and Sons,Inc.,1989
\bibitem{ES} A.Enge and R.Schertz, 
 Constructing elliptic curves over finite fields using double eta-quotients, 
J.Th{\'e}or.Nombres Bordeaux 16 (2004), 555--568.
\bibitem{GA} A.Gee, Class invariants by Shimura's reciprocity law, J.Th{\'e}or.Nombres Bordeaux 11 (1990),45-72.
\bibitem{II} N.Ishida and N.Ishii, Generators and defining equation of the modular function field of the group $\Gamma_1(N)$, Acta Arith. 101.4 (2002),303-320.
\bibitem{I1} N.Ishii, Rational expression for $J$-invariant function in terms of generators of modular function fields, Int.Math. Forum 2 (2007) no. 38, 1877 - 1894.
\bibitem{LA}
S.Lang, Elliptic Functions, Springer-verlag,1987.
\bibitem{RS}
R.Schertz,
Weber's class invariants revisited, 
J.Th{\'e}or. Nombres de Bordeaux 14(1) (2002), 325--343.
\bibitem{SG}
G.Shimura, Introduction to the Arithmetic Theory of Automorphic Functions, Iwanami-Shoten and Princeton University Press,1971.
\bibitem{SJ} J.Silverman, The Arithmetic of Elliptic curves, Springer-verlag,1986.
\bibitem{YAI}
S.Yoshimura,A.Comuta and N.Ishii, $N$-systems,class polynomials of double eta-quotients and singular values of $j$-invariant function, Int.Math.Forum 4 (2009) no.8,367-376.
\end{thebibliography}
\end{document}